\documentclass[12pt]{article}
\usepackage{amssymb}
\usepackage{epsf}

 \newcommand{\EE}{\mathbb E}
 
 \newcommand{\RR}{\mathbb R}
 \newcommand{\CC}{\mathbb C}
 \newcommand{\ZZ}{\mathbb Z}
 \newcommand{\II}{\mathbb I}
 \newcommand{\TT}{\mathbb T}
 \newcommand{\QQ}{\mathbb Q}
 \newcommand{\HH}{\mathbb H}

\newcommand{\vol}{{\rm vol}}

\newcommand{\h}{\mbox{\small $\#$}}

\newtheorem{theorem}{Theorem}

 \newcommand{\be}{\begin{equation}}
 \newcommand{\ee}{\end{equation}}

\begin{document}

\title{Model Sets: A Survey}

\author{Robert V. Moody}

\maketitle

\begin{abstract}\par
This article surveys the mathematics of the cut and project
method as applied to point sets, called here {\em model sets}.
It covers the geometric, arithmetic, and analytical sides
of this theory as well as diffraction and the connection with 
dynamical systems.
       \end{abstract}

%%\begin{minipage}{0.8\textwidth}
\bigskip\bigskip
\begin{quotation}

\footnotesize{
\centerline{Dedicated to the memory of Richard (Dick) Slansky}
\bigskip
The spirit of the universe is subtle and informs all life.
Things live and die and change their forms, without knowing
the root from which they come. Abundantly it multiplies;
eternally it stands by itself. The greatest reaches of space
do not leave its confines, and the smallest down of a bird
in autumn awaits its power to assume form. \hfill  --- Chuang Tzu 
(tr. Lin Yutang)}

\end{quotation}
%%\end{minipage}

\section{Introduction}

Even when reduced to its simplest form, namely that
of point sets in euclidean space, the phenomenon of 
genuine quasi-periodicity appears extraordinary.
Although it seems unfruitful to try and define the concept
precisely, the following properties may be 
considered as representative:
\begin{itemize}
\item discreteness
\item extensiveness
\item finiteness of local complexity
\item repetitivity
\item diffractivity
\item aperiodicity
\item existence of exotic symmetry (optional).
\end{itemize}

The purpose of this paper is to give an overview of
the mathematics of the cut and project method which 
not only provides
a very rich harvest of point sets (called model
sets) satisfying these properties, but
also provides a very natural way to link these ideas with many other
structures in mathematics.

A subset $\Lambda$ of $\RR^d$ is called a {\em Delone (Delaunay)} set
if it is uniformly discrete and relatively dense. This means that
there are radii $r,R > 0$ so that each ball of radius $r$ (resp. $R$)
contains at most (resp. at least) one point of $\Lambda$. Although
this is a fairly strong version of the first two items on our list,
it is the most commonly used one and coincides well with the 
primitive atomic picture of a (ideally infinite)
piece of material. 

The set $\Lambda$ has {\em finite local complexity} if
for each $r >0$ there are, up to translation, 
only finitely many point sets (called {\em patches
of radius} $r$)
of the form $\Lambda \cap B_r(v)$.
Here $B_r(v)$ is the ball of radius $r$ about the point $v \in \RR^d$.
So, on each scale, there are only finitely  many different patterns
of points. This condition can be expressed topologically:
$\Lambda$ has finite local complexity iff the closure
of  $\Lambda -\Lambda$ is discrete.
It is conceivable to replace ``translation'' by ``isometry''
in this definition, but the theory would change considerably and,
with the notable exception of the pinwheel tiling \cite{radin},
little has been said so far on this more general situation.

Repetitivity means loosely that any finite patch that appears,
appears infinitely often. More precisely, given any
patch of radius $r$ there is
an $R$ so that within each ball of radius $R$, no matter its 
position in $\RR^d$, there is at least one translate of this patch. 
A stronger form of this requires in addition that each type
of patch of radius $r$ should appear with a well-defined
frequency. The sets that we deal with here normally have
this additional property (see Section 3).

From the very beginning, diffractivity has been the hallmark
of aperiodic order. Physically it is the most visible of its manifestations.
Mathematically it is one of the most subtle and least visible! Very roughly
we are asking (mathematically) that the Fourier transform of the 
autocorrelation density that arises by placing a delta peak
on each point of $\Lambda$, should contain a part that looks
discrete and point-like. Later on we will make this very vague
prescription precise. The amazing thing is that in the context
of model sets we obtain perfect
diffractiveness, in the sense that the diffraction is purely point-like,
under the fairly mild  hypothesis of regularity. One of the 
goals of this survey is to show how this comes about.

Lack of periodicity speaks for itself. Lattices and unions of
cosets of lattices are the basis of the most prevalent forms
of long-range order (crystallography). Point sets based on them 
satisfy all the previous properites. But of course, that is the trivial part
of the theory! The objective is to move into new territory.

Exotic symmetry usually means non-crystallographic symmetry.
Although not mathematically essential, certainly the existence of physical
structures with ``forbidden'' icosahedral symmetry was instrumental
in the rapid development of this field.

Even with the rather strong interpretations on the 
various properties listed above, we still do not
know how to characterize sets that satisfy them. For 
an extensive discussion of these problems see \cite{Lagarias2}.
However there is one very general method of construction 
which relies on controlled projection from a discrete group located in 
some auxilliary ``embedding'' space. 
In its original form this so-called cut and project method is based on 
projection from lattices in higher dimensional spaces. Many people
have written about this starting, in physics, with
the work of P.~Kramer \cite{Kramer} and including the very useful 
article of
Y.~Meyer given in the previous edition of this School 
\cite{YM2}.
Meyer had already thought about sets formed by projection
from the view point of harmonic analysis long before the 
discovery of quasi-crystals \cite{YM1}. Even though it is convenient
to somewhat rearrange the main components of his original construction,
nonetheless he created a formalism
which is ideal for creation of points sets with the desired
properties of long-range aperiodic order. These are the 
{\em model sets}. 

Some people object to the terminolgy  ``model set''
prefering ``cut and project set'' which sounds more serious
and professional. However,
we prefer to interpret ``model'' as meaning exemplary and think that
in terms of both of its priority and its greater generality the term 
deserves to be adopted. 

The main purpose of this article is first to give some 
idea of the scope of the relevant examples that arise
as model sets (this scope surely not yet fully realized)
and then to show how the model sets are poised between 
a number of quite different areas of mathematics. It is the satisfying
way in which they connect many diverse parts of mathematics
that makes model sets so intriguing and offers to the imagination
so many tantalizing prospects for future work. For the reader interested
in more on the tiling side of quasiperiodicity we recommend the survey
paper \cite{michaelSurvey}, which also provides a complementary source of some
of the material presented here.

\section{Model sets}

Let us launch ourselves directly into the notion of a model set. 
By definition, a cut and project scheme 
consists of a collection of spaces and mappings

\be \label{cutandproject}
  \begin{array}{ccccc}
   \RR^d & \stackrel{\pi^{}_1}{\longleftarrow} & \RR^d \times G &
           \stackrel{\pi^{}_2}{\longrightarrow} & G  \\
    & & \cup & & \\ & & {\tilde L} & & \end{array}
\ee
where $\RR^d$ is a real euclidean space and $G$ is some
locally compact abelian group, $\pi^{}_1$ and 
$\pi^{}_2$
are the projection maps onto them, and ${\tilde L} \subset 
\RR^d \times G$
is a lattice, i.e., a discrete subgroup such that the quotient group
 $(\RR^d\times G)/{\tilde L}$ is compact. We assume that
 $\pi^{}_1|^{}_{\tilde L}$ is
 injective and
that $\pi^{}_2({\tilde L})$ is dense in $G$. We call $\RR^d$
(resp.\ $G$) the {\em physical} (resp.\ {\em internal}) space.
The product $\RR^d \times G$ is the {\em embedding} space. We 
write $ L = \pi^{}_1({\tilde L})$. It is very convenient
to define the mapping
\be
 ^* : L \longrightarrow G \; : \; \mapsto \pi_2(\pi_1|_L^{-1}(x)) \;.
\ee

Given any subset $W \subset G$, we define a corresponding 
set $\Lambda(W) \subset \RR^d$ by
\be
  \Lambda (W)\; = \; \{\pi_1(x)  \;|\; x \in {\tilde L} , \, \pi_2(x) \in W \}
 \; = \; \{u \in L \;|\; u^* \in W \} \,.
\ee 
We call such a set $\Lambda$ (or more generally any 
translate of such a set) a {\em model set} 
(or {\em cut and project set}) if the following condition
[W1] is fulfilled:
\begin{itemize}
\item[\bf W1:] $W$ is nonempty and 
$W \;  =  \; \overline{\mbox{int}(W)} \;$ is compact.
\end{itemize}

For some of the deeper results we need more precise assumptions,
of which the following are the most relevant:
\begin{itemize}
\item[\bf W2:] The model set $\Lambda$ is {\em generic} if
the boundary $\partial W$ of its window $W$ 
satisfies is $\partial W \cap \pi_2({\tilde L}) = \emptyset$.
\item[\bf W3:] The model set 
 $\Lambda$ is {\em regular} \footnote{The terminology here is not
standardized. Sometimes what we call generic is called regular.
Nonetheless,
the generic situation is justifiably ``generic'' as we point out below. 
Our use of regular is close to the one used in \cite{martin2,martin3}.}  
if $\partial W$ is of (Haar) measure $0$.
\end{itemize}

The definition formalizes the notion of a point set
in $\RR^d$ constructed by projecting selected points
from a lattice in some ``super--space''. The points
selected for projection are those which fall into some
bounded region when they are projected into the complementary
internal space $G$.  The notion of lattice, familiar in a real
space as the $\ZZ$-span of a basis of that space, is replaced
here by the more general definition that can be applied
to any topological group. In condition [{\bf W1}] the equality
could be replaced by ``$\subset$'', but it is convenient
to have this additional hypothesis since then $\overline{\Lambda^*}
= W$.

We asssume that the reader is familiar with the 
the most common, and only easily visualized examples of this,
that are based
on a pair of orthogonal axes, at irrational slopes (i.e.
axes, through the origin of the standard lattice $\ZZ^2$ in $\RR^2$
which are taken to be the physical and internal spaces, 
and a window which is an interval on the internal space
e.g. \cite{Senechal, michaelSurvey}.

There are four different view points to the diagram above,
which we can picture as follows. In the first place we have
the physical space $\RR^d$ and the point set $\Lambda$ in
it whose geometric properties are those that we wish to
understand and describe. Lattices are discrete groups
inside larger continuous groups, and so may be thought
of as {\em arithmetic} in origin. We will see later how in many
interesting cases the arithmetic aspect is quite central.
Why the internal side should be thought of as having to do
with {\em analysis} will also emerge in later, but an initial
way to think of it is that on the internal side, the set
of points of $\Lambda$ appear in a totally different arrangement
so that their closure is a very nice region of space.
Finally, by definition, $\TT := (\RR^d \times G)/{\tilde L}$ is a
compact abelian group. In the usual situation of a real internal space,
$\TT$ is a torus, whence the notation. In any case,
$\TT$ has a totally natural action of $\RR^d$ on it and it is
this action that gives rise to a {\em dynamical system}. In the end
we will have a second, related, dynamical system which plays an important role
in questions around diffraction.

\be \label{viewpoints}
  \begin{array}{ccccc}
& & \mathrm{dynamical \ systems \ side} & & \\
 & & \TT &  & \\
 & &  \uparrow & &\\
   \RR^d & \longleftarrow & \RR^d \times G &
           \longrightarrow & G  \\
  \mathrm{physical \ side}  & & \uparrow & & \mathrm{analytical \
 side}\\
 & &\tilde{L} & & \\
 & & \mathrm{arithmetic \ side} & & \\
\end{array}
\ee

In fact this picture can be dualized, thereby producing
yet another four pictures! This dualization plays quite
an important role in Meyer's theory which we touch on
only most briefly here. However we will use one part of the
dual picture. In the dual picture it is ${\hat \TT}$ that 
is the lattice and we have 
\be \label{dualProj}
  \begin{array}{ccccc}
  \widehat{\RR^d} & \stackrel{\hat{\pi}_1}{\longleftarrow} & 
\widehat{\RR^d} \times \widehat{G} &
           \stackrel{\hat{\pi}_2}{\longrightarrow} & \widehat{G}  \\
    & & \uparrow & & \\ & & \hat{\TT} & & \end{array}
\ee
Here we have identified the dual of the direct product $\RR^d \times G$
with the direct product of the duals, and have chosen to single
out the canonical projections as the important maps, designating
them by $\hat{\pi}_1$ and  $\hat{\pi}_2$ respectively.

\section{Geometric side} 
The geometric properties of model sets $\Lambda = \Lambda(W)$ 
have been described
in detail elsewhere (for instance \cite{RVM, martin2}). 
We will restrict ourselves to pointing out a few of the
most important features here. In the first place, model
sets are Delone sets and have the property of finite local
complexity. In fact they satisfy a very strong form
of finite local complexity:

\be \label{meyer}
\Lambda - \Lambda  \; \mbox{is uniformly discrete}.
\ee

A Delone set satisfying (\ref{meyer}) is called a {\em Meyer
set}. There are a remarkable number of ways describing
Meyer sets (\cite{YM2, RVM, Lagarias1}) which link them strongly
with harmonic analysis. 
Though the Meyer property is considerably  weaker than
that of a model set, we nonetheless have

\begin{theorem}(\cite{YM1})\label{MeyerTheorem}
Any Meyer set is a Delone subset of some model set.
\end{theorem}

The situation regarding repetitivity is
complicated by the boundary of the window $W$. 
$\Lambda$ is repetitive if it is generic. If we are allowed
to modify a model set by moving its window around then
it is straightforward using the fact that
$G$ is a Baire space to see that the window can
be moved to make the resulting model set generic.
A proof of this can be found in \cite{padics}.
Furthermore, in the regular case, the frequency of repetition
 of each patch
is well-defined in the sense that for each finite patch
$P$ the number of occurences
of the patch $P$ (up to translation) per unit of volume in the ball
$B_r(0)$ of radius $r$ approaches a positive limit
as $r \to \infty$. This is actually not hard to prove
once one has established uniformity of projection (Theorem \ref{uniform}). 

Lack of periodicity is automatic for model sets as
long as the mapping $^*$ is injective. Otherwise, 
the kernel of $^*$ is the translation group of $\Lambda$.

The comprehensive paper of Lagarias \cite{Lagarias2} is the most
extensive study to date of the geometry of point sets
in the context of quasi-periodic structures.

\section{Arithmetic side}

Although the requirement in the definition 
of a model set of the existence of a lattice is not
in itself particularly arithmetic, nonetheless the 
interesting and important examples all have strong
arithmetic aspects. In the usual cases where the internal
space is a real space, the arithmetic arises through 
the standard inner product on the embedding space and the
nature of the two projections. 

We will illustrate here the typical arithmetic input into
the theory with two very different examples. 

\subsection{The icosian model sets}

It was M.~Baake et al \cite{michael} who first pointed out that 
the root and weight
lattices of types $A_4$ and $D_6$ could be used
as the lattices for projection in cut and project schemes
for dihedral ${\cal D}_5$ and icosahedral symmetries in $2$ and $3$
dimensions. The fact that these two groups are Coxeter
groups (finite reflection groups) and form the first two of the
the series $H_2,H_3,H_4$ of non-crystallographic finite Coxeter
groups \footnote{The group $H_2$ is the dihedral group of
order $10$. Usually it is fitted into the series $I_2(k)$
of dihedral Coxeter groups, but it is also completely natural
to think of it in the icosahedral series, as we do here. In fact,
we could go a step further and include $H_1$ which is simply the
reflection group of order $2$.} 
(of which $H_3$ and $H_4$ are the only examples of rank
larger than $2$) suggests that $H_4$ should also appear
in this context. This was first pointed out in \cite {ES}
and elaborated in more detail in \cite{MP}. We do nothing
more than outline this here. It is not necessary to
know anything about root systems to follow this example.

The elements of norm $1$ of the usual quaternion ring  
$\HH = \RR + \RR i +\RR j + \RR k$ form a group isomorphic
to ${\rm SU}(2)$. Since this group is a $2$-fold cover of
the orthogonal  group ${\rm SO}(3)$, in particular it contains $2$-fold
covers of the icosahedral group.  One such example is
the following list $I$ of $120$ vectors:

\be
\begin{array}{c}
\frac12(\pm1,\pm1,\pm1,\pm1),\quad (\pm1,0,0,0)\   
                                  \textrm{and all permutations},\\
\frac12(0,\pm1,\pm\tau',\pm\tau) \quad 
                                  \textrm{and all even permutations}
\end{array} 
\ee
where $\tau = (1+\sqrt5)/2$ is the Golden ratio and $'$ indicates
the conjugation map $\sqrt{5} \mapsto -\sqrt{5}$.

The subring $\II$ generated by this group 
is called the {\em icosian ring}. Of course it depends
on our particular choice of $I$, though it is straightforward
to see
that $\II$ is unique up to inner automorphisms of $\HH$.
The form of the points of $I$ makes it clear that
$\II$ is a $\ZZ[\tau]$-module.
We let $^*$ denote the mapping
on $\II$ that conjugates each of the coordinates
with respect to the unique Galois non-trivial
automorphism on $\ZZ[\tau]$ (defined by sending 
$\sqrt5 \mapsto -\sqrt5$). Note that $\II^* \ne \II$.

The ring $\II$ is of rank $4$ over $\ZZ[\tau]$ and
rank $8$ over $\ZZ$. We make an explicit embedding
of $\II$ as a lattice ${\tilde \II}$ in $\RR^8$ by the 
mapping $x \mapsto (x,x^*)$.

This already provides the framework
of a cut and project scheme:
\be 
  \begin{array}{ccccc}
   \RR^4 & \stackrel{\pi^{}_1}{\longleftarrow} & \RR^4 \times \RR^4 &
           \stackrel{\pi^{}_2}{\longrightarrow} & \RR^4  \\
    & & \cup & & \\ & & {\tilde \II} & & \end{array}
\ee
with the projections being given by the first and second
components of $(x,x^*)$.

Remarkably the lattice ${\tilde \II}$ has an entirely
natural interpretation as the root lattice of type $E_8$
(see for instance \cite{ConwaySloane}
which underscores its arithmetic nature. This is explained
in \cite{ES,MP,MP2}.   

Now we wish to show that this cut and project
scheme respects the symmetry that is inherent in its
construction.  Geometrically the points of $I$ form the vertices 
of a regular
polytope $P$ in $4$-space and also form the vectors of a root
system $\Delta_4$ of type $H_4$. 
The Coxeter group $H_4$ is none other than the group of automorphisms
of $P$ (and also of $\Delta_4$), and is in fact
very easily described: it is the set of all ($14400$) maps
\be
 x \mapsto uxv \quad ; \quad  x \mapsto u{\overline x}v
\ee 
where $u,v \in I$.
The subgroup of these transformations in which $v= u^{-1}$
is obviously a copy of the icosahedral group and 
this subgroup stabilizes the 3-dimensional space
$\RR i + \RR j + \RR k$
of {\em pure} quaternions.  

These maps provide automorphisms of the rings $\II$ and,
via conjugation, on  
$\II^*$ too, and thus give rise to an action of $I$ as
automorphisms on the
entire cut and project scheme. If the window $W$ is chosen
to be invariant under $I$ then the resulting model
set is also $I$-invariant. 
 
Restricting everything to the pure quaternions we get a new 
cut and project scheme based on the $6$-dimensional root lattice
$D_6$ and an icosahedral
symmetry. Restricting further to the 
planes orthogonal to the 5-fold axes brings us 
back to $A_4$ and the related dihedral ${\cal D}_5$ symmetry. 
A step further, and we arrive at the Fibonacci chain in 
$1$ dimension.
Thus all three families as well as the fundamental Fibonacci
model sets fit together in this quaternionic model. Not only
is this very pretty, it also essentially 
encompasses the generic
situation for icosahedral symmetry in model sets:
the only other relevant lattices in $6$-space
are the $D_6$ weight lattice and the lattices
lying between the root and weight lattices. For more
on this see \cite{Merminetal,MP2}.

\subsection{$p$-adic model sets}

Until recently, little thought had been given to the situation
in which the internal group is something different than
another real space, or at worst a real space crossed
with a torus.  However, there is
a whole series of very natural locally compact
abelian groups that are not euclidean in nature, namely
the $p$-adic groups. Since these may not be familiar
in this context let us recall the basic ideas. 

Let $p$ be a prime number in the integers $\ZZ$. Using $p$ we can 
define a metric on the rational numbers $\QQ\/$, and
by restriction on $\ZZ$, in the following way. For each $a \in \ZZ$,
we define its $p$-value, $\nu_p(a)$, as
the largest exponent $k$ for which $p^k$ divides $a$ 
(with $\nu_p(0):= \infty$). This function
is extended to the $p$-adic valuation $\nu_p: \QQ \longrightarrow \ZZ$ by 
$\nu_p(a/b) := \nu_p(a) - \nu_p(b)$ for all rational numbers
$a/b$. We now define the ``distance'' between two rational numbers $x,y$
as $ d(x,y) = p^{-\nu_p(y-x)}$ .

It is not hard to see that this does define a metric on $\QQ$, in which 
closeness to $0$ is equivalent to high divisibility by the prime
$p$. The completion of the rationals under this topology
is the field of $p$-adic numbers ${\QQ_p}$ and the 
completion of the subring $\ZZ$ is the subring of $p$-adic 
integers, ${\ZZ_p}$.
Each $p$-adic integer can be given the more concrete
representation as a series
in the form $\sum_{n=0}^\infty a_np^n$
where the $a_n$ are integers in the range $0 \le a_n < p$.
Note that convergence here is automatic, even though there
are infinitely many terms in the sum, because of the nature
of the $p$-adic topology.
The topologies defined by such metrics 
have other counter-intuitive properties. For example,
for each non-negative integer $k$, the set $p^k\cdot {\ZZ_p}$, 
of elements of ${\ZZ_p}$ divisible by $p^k$, is the ball 
of radius $p^{-k}$ and is clopen, i.e.\ both open and closed,
as too are all its cosets, $a + p^k\cdot {\ZZ_p}$.
 
Seen as a topological space, ${\ZZ_p}$
is both compact and totally disconnected
(but not discrete).
In particular, ${\QQ_p}$ and ${\ZZ_p}$
are locally compact abelian groups under addition.
Thus, we can use ${\ZZ_p}$ 
to construct interesting cut and project 
schemes for $\RR^d$ simply by taking $G := ({\ZZ_p})^d$ and $L = \ZZ^d$  
embedded diagonally into $\RR^d \times {\ZZ_p}^d$
(based on the natural embedding of $\ZZ$ in $\ZZ_p$).
For more on $p$-adic numbers and other totally disconnected
groups, the reader may consult \cite{Neukirch, Bourbaki}.

In \cite{padics} it was shown that 
a number of interesting substitution systems and tilings
can be interpreted in a $p$-adic setting, including the well-known
chair tiling. Rather than repeat these, let us give a different example,
mentioned in \cite{padics} but not elaborated upon. 
One of the
earliest classes of aperiodic tilings to be discovered was
the class of Raphael Robinson's
square tilings \cite{Robinson}. The title of his paper
recalls that the mathematical interest in aperiodic structures
had a totally different (and earlier) origin than the physical one, namely
the interest in decidability problems in the tiling of the plane
with tiles of finitely many different types.
The Robinson tilings are tilings of the plane by equally sized squares, in 
the usual fashion, with the twist that the square tiles come in 6 types 
(up to rotational and reflectional symmetries), distinguished by the 
markings of their edges,
and the tiling is required to respect these edge markings by having
the edges of adjacent tiles properly matched.
Pictures of the tiles may be found in \cite{Robinson} and \cite{GS}.
What is important for our discussion is that there is another set of 
markings
by lines of these tiles and the correct tilings are those for which
these lines arrange themselves into a pattern of squares of increasing scales 
$1,2,4,8, \dots$ (see Fig.1).  This picture is the one that Robinson 
used to prove the
aperiodicity, for evidently no translation can map the squares of all scales
onto themselves simultaneously. The same idea was used by Penrose in
his recent hexagonal tiling \cite{Penrose}. 
\vspace*{8mm}
\begin{figure}[ht]
\centerline{\epsfysize=65mm \epsfbox{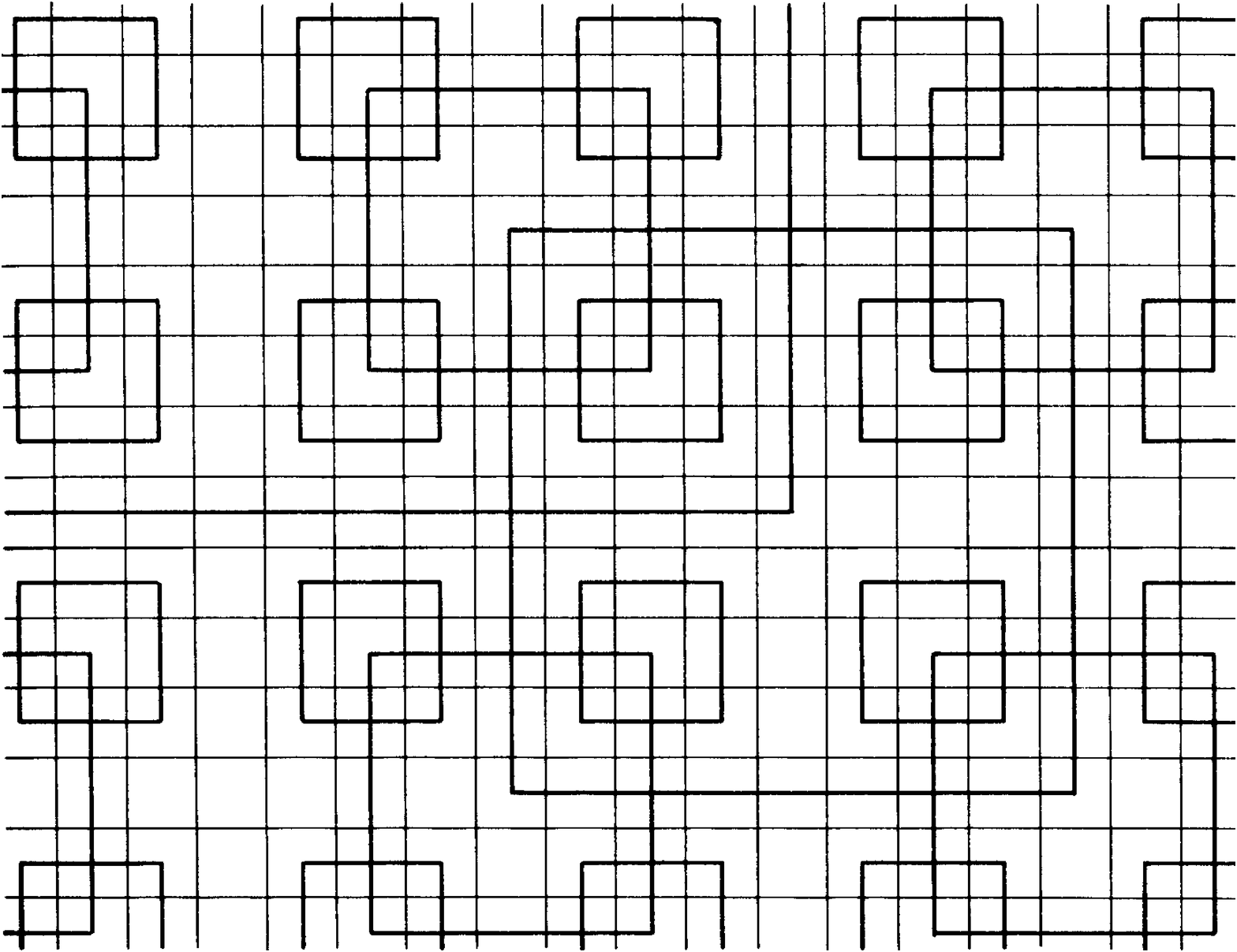}}
\caption{Pattern of increasing squares in a Robinson square
tiling}
\end{figure}

Now the point is that the centres of the tiles of each of the six
types form a model set based on an internal space which is $2$-adic. 
Very briefly the argument is as follows.

Starting with Fig.~1 and {\em ignoring the actual square tiles themselves}
we have patterns of interlocking squares of increasing scale.
Let's call these the {\em pattern-squares} to distinguish them
from the actual tiling squares.
The vertices of the smallest scale pattern-squares (order 1) form the 
vertices
of a lattice $L'$ once one of them, say $c_1 =(0,0)$, has been chosen
as the origin. It is convenient to introduce the larger
lattice $L = \frac12 \cdot L'$ which we can identify with $\ZZ^2$.
The locations of the vertices of 
the increasingly scaled pattern-squares 
determine two sequences,
$\{\alpha_k\}$ and   $\{\beta_k\}$, composed out
of the two numbers $\{\pm1\}$ as follows:
the vertices of the pattern-squares
of order $k$ (side-length $2^k$) are the points of a coset
\be \label{vertices}
L_k := c_k + 2^k L.
\ee 

where 
\be \label{sequence}
c_k = (\alpha_0 +\alpha_1 2 + \dots , + \alpha_{k-2} 2^{k-2},
\beta_0 +\beta_1 2 + \dots , + \beta_{k-2} 2^{k-2})
\ee
for $k = 2,3,  \dots$. Conversely, given two sequences $\alpha,
\beta$ of $\pm 1$'s we can use (\ref{sequence}) and (\ref{vertices})
to define the vertices of a suitable pattern of squares.

The pattern-squares themselves are determined by the condition
that the points of $L_k$ are the centres of the squares
of order $k-1$ for all $k >1$. This then establishes a coordinatization
of the pattern-squares. Both the actual pattern of the
pattern-squares
and their coordinatization depend on the choice of
our two sequences (more below).

Now looking again at the tiling squares,
we distinguish the $6$ types of tiles according to
their location in the pattern-squares. We list these here
together with a description of the coordinates of the centres of
their tiles:
\begin{itemize}
\item[(1)] the ``corner tiles'' of the squares of order $1$;
\newline 
coordinates  $L_1 = 2L$;~ density $\frac14$;
\item[(2)] the ``corner tiles'' of all squares of all higher orders;
\newline
coordinates $\bigcup_{k\ge 2} c_k + 2^kL$;~
 density $\frac{1}{12}$
\item[(3)] ``cross tiles'' where edges of two different orders of pattern
squares meet (actually these orders always differ by exactly $1$); 
\newline coordinates:
\newline 
$\left( \bigcup_{k \ge 3}c_k + (\pm 2^{k-2},\pm 2^{k-3})+ 2^kL \right)
\cup 
\left( \bigcup_{k \ge 3} c_k+ (\pm 2^{k-3}, \pm 2^{k-2}) + 2^kL \right)$;
\newline 
density $\frac16$; 
\item[(4)] ``edge '' squares, which contain part of a single edge of a 
pattern-square, except those in which are exactly in the middle of an edge;~
density $\frac16$;
\item[(5)] ``edge '' squares, which contain part of a single edge of a 
pattern-square and which are exactly in the middle of an edge;~ 
density $\frac16$;
\item[(6)] blank tiles, with no part of any edge in them;~ 
density $\frac16$.
\end{itemize} 

Observe that each of these sets is a countable union $w_j$ of cosets of 
the form $a + 2^k L$. Let us replace each of these by the corresponding 
$2$-adic clopen set $a + 2^k\cdot \ZZ_2^2$. In this way
we get $6$ open sets $W_j,\; j = 1,\dots ,6$, whose closures, being
closed subsets of the compact group ${\ZZ_2^2}$, are compact
with non-empty interiors.
Finally we can describe
the centres of the squares of type $j$ as
\be
\{ x \in \ZZ^2 \mid x \in W_j \}
\ee
which, unlikely as it appears, is a model set under the scheme

\be \label{padiccutandproject}
  \begin{array}{ccccc}
   \RR^2 & \stackrel{\pi^{}_1}{\longleftarrow} & 
\RR^2 \times {\ZZ_2^2} &
           \stackrel{\pi^{}_2}{\longrightarrow} & G  \\
    & & \cup & & \\ & & \ZZ^2 & & \end{array}
\ee
where $\ZZ^2$ is embedded into $\RR^2 \times {\ZZ_2^2}$
diagonally: $x \mapsto (x,x)$.

The entire tiling is determined by the vertices of the 
various squares and hence by the window
\be
W := \bigcup_{k=1}^\infty c_k + 2^k\ZZ_2^2.
\ee
Evidently $W$ is an open subset of ~ $\ZZ_2^2$. Let 
$w \in \overline{W}\backslash W$. Then for each $m \in \ZZ_+$,
$w+2^m ~ \ZZ_2^2$ meets $W$, so $w \equiv c_k \; {\rm mod}\;
2^{{\rm min}\{m,k\}}~ \ZZ_2^2$ for some $k = k(m)$.
If $k \le m$ then $w \in c_k + 2^k\ZZ_2^2 \subset W$.
Thus $k > m$ and $w \equiv c_{k(m)}{\rm mod}\;2^m \;\ZZ_2^2$.
It follows that $w$ is the limit of some subsequence of
the $\{c_k\}$. Since the entire sequence evidently converges
(in the $p$-adic topology, of course!) to some $c = (a,b)$,
where $a = \sum \alpha_k 2^k$ and similarly  for $b$, 
we see that $w= c$ and so $\partial W = \{c\}$. Thus
the model set of all vertices is regular, and even generic
provided that $c \notin \ZZ^2$.

More generally, one may expect this $p$-adic topologies to arise
whenever there is a self-similarity 
$\theta : L \rightarrow L$
for which $\theta(L) \subset L$, but $\theta(L)\neq L$. 

In \cite{padics} we also see the appearance of mixed $p$-adic and
real spaces as the internal spaces. Beyond these types we are not
aware of any interesting examples, though they may well exist.

\section{Analytic side}

The transition from an inherently discete picture
on the physical side to something inherently
far more continuous on the
internal side is made via H.~Weyl's theory of uniform
distribution. 

Let us assume that we have a model set $\Lambda$.
 Now consider the following question. Suppose
that we take a ball $B_R(0)$ of radius $R$ about the
origin in $\RR^d$ and look at $\Lambda_R := \Lambda \cap B_R(0)$.
Then we can ask how $\Lambda_R^*$ is distributed over $W$.
We say that the sets  $\Lambda_R^*$ are {\em uniformly distributed} 
if for each open set $U \subset W$ we have 
\be
\lim_{R \to \infty} \frac{{\rm card}(\Lambda_R^* \cap U)}{\mu(W)}
\; = \; \mu(U)/\mu(W)
\ee
where $\mu$ is Haar measure on $G$. 

\begin{theorem}\cite{martin2,bert3} \label{uniform}
If $\Lambda$ is regular then the sets $\Lambda_R^*$
are uniformly distributed over $W$.
\end{theorem}

Let $f^*:G \longrightarrow \CC$ be any function.
We define
$f:L \longrightarrow \CC$ 
by $f(x) = f^*(x^*)$. If $f^*$ is supported on the window $W$
then evidently $f$ is supported on the
model set $\Lambda$. If $f^*$ is continuous (which is the
case of interest) then this is iff. 

\begin{theorem} (Weyl) \cite{weyl}
 If $\Lambda$ is regular and  $f^*$ is continuous
then
\be 
\lim_{R\to \infty} \frac{1}{{\rm card}(\Lambda_R)}
\sum_{x \in \Lambda_R} f(x) 
\; = \; \frac{1}{{\rm vol}(W)}\int_{W} f^*(u) d\mu(u)
\ee
\end{theorem}

Since $W$ has boundary of measure zero, it is not
necessary to insist that $f^*$ (which is supported on $W$)
be continuous
on all of internal space, only on the window $W$.

In this way discrete averaging on the model set is transformed into
integration on the window. This process was used in \cite{BM}
to determine the existence of invariant measures on internal
space in the presence of self-similarity on the quasi-crystal.
We briefly explain this. We assume here that internal space
is $\RR^n$ for some $n$.

A {\em self-similarity}
of $\Lambda$ is an affine linear mapping $t=t^{}_{Q,v}$
\be \label{selfsim}
     t^{}_{Q,v} \; : \quad x \mapsto Qx + v
\ee
on $\RR^d$ that maps $\Lambda$ into itself, where $Q$ is a (linear) 
similarity
and $v\in\RR^d$. Thus $Q = q R$, i.e.\ it is made up of an orthogonal 
transformation $R$ and an {\em inflation factor} $q$.

Let $t^{}_{Q,v}$ be a self-similarity of $\Lambda$. Since $\Lambda$ is
uniformly discrete, we must have $|q|\geq 1$. We will assume
$|q|>1$ and that $QL=L$.
We are interested in the {\em entire} set
of affine inflations with the same similarity factor $Q$.

Note that $Q$ naturally gives rise to an automorphism $\tilde{Q}$ of
the lattice $\tilde{L}$, i.e. an element of ${\rm GL}_{\ZZ}(\tilde{L})$,
and a linear mapping $Q^*$ of $\RR^n$ that maps
$W$ into itself. From the arithmetic nature of $\tilde{Q}$
we deduce that the eigenvalues of $Q$ and $Q^*$ are algebraic
integers and from the compactness of $W$ that 
$Q^*$ is contractive. 

Define 
\be 
      W^{}_Q \; := \; \{ u \in \RR^n \mid
                Q^* W + u \subset W \} \, ,
\ee
We say that $Q$ is {\em compatible} with $\Lambda$ if
$\mbox{int}(W^{}_Q)\neq\emptyset$. Assuming that this
is the case (not a strong assumption) , 
then the set ${\cal T}_Q$ of
affine inflations with the same similarity $Q$ is the set of mappings
$t^{}_{Q,v} : x \mapsto Qx + v$, where $v$ runs through the set
\be 
     T \;= \; T^{}_{Q} \; := \; \{ v\in L \mid v^* \in W^{}_Q \} \, .
\ee

\begin{theorem}
If $Q$ is a self-similarity
and the above assumptions on $Q$ apply then there is a unique 
absolutely continuous positive measure $\mu$ on internal space,
supported on $W$, satisfying:
\begin{itemize}
\item{}$\mu^{}(W)= 1$;
\item{}
$\mu^{}$ 
is invariant in the sense that, if we define $t^*_v \cdot\mu^{}_f$
by  $t^*_v \cdot\mu(Y) = \mu((t^*_v)^{-1}(Y))$, then
\be \label{invariantdensity}
   \mu^{} \; = \; \lim_{s \to \infty} \frac{1}{\h \left(T \cap B_s(0)\right)}
\sum_{v\in T \cap B_s(0)} t^*_v\cdot\mu^{}_f ~.
\ee
\end{itemize}
\end{theorem}
 The similarity of this measure
to Hutchison measures in the context of iterated function systems is not
coincidental. In fact, if we restrict to the ball $B_s(0)$
then the $\{ t_v^*\}$ form a finite set of contractions
which is indeed an iterated function system.
For more on this and invariant density functions 
on model sets see \cite{BM}. 

The type of limit averaging involved here is a very natural one
from the point of view of physical situations, representing
the transition from the world of sets finite in extent to the ideal 
world of infinitely extended point sets. 

Although no one to our knowledge has made any use of it,
it is interesting to use Weyl's theorem to transfer the
structure of $L^2(W)$ to a space of similar objects
on $\Lambda$. Namely, the space of continuous functions
on $W$ leads to a space 
\be \label{pseudoContinuous} 
{\cal C} = {\cal C}(\Lambda)
\ee
 of funtions 
on $\Lambda$ via
the mapping $^*$. Then the usual inner product
$\langle f, g\rangle = \int_W \overline{f^*(u)} g^*(u) du 
= \langle {\overline f}, g\rangle_W$ defines an inner 
product on ${\cal C}$ and we can complete this space in
order to get a Hilbert space $\overline{\cal C}$ isomorphic 
to $L^2(W)$. Of
course the elements of $\overline{\cal C}$ can no longer be
interpreted as functions on $\Lambda$ since
functions on $\Lambda$ that differ by a function whose
absolute square has limit average sum equal to $0$
are identified. 

\section{Dynamical systems side}

So far we have looked at one model set in isolation.
Now we move on to consider families of model sets. We start
with a number of definitions and results. 
All of these may be found in the paper of
M.~Schlottmann \cite{martin3} on which
we have relied heavily here. Many are well-known
in the context of tilings for which a recent reference with a good
bibliography is \cite{boris}. In this section all point
sets under discussion are assumed to be Delone sets in $\RR^d$.

Two Delone sets $S,S'$ in $\RR^d$ are {\em locally isomorphic}
(or some people say {\em locally indistinguishable})
if, up to translations, every patch of either of them occurs
in the other. Thus on any finite scale, up to translation, the two sets
are indistinguishable. Given a 
Delone set $S \subset \RR^d$ we can
look at its local isomorphism class (LI class) ${\rm LI}(S)$, 
namely all point sets locally isomorphic to it.

We denote by ${\cal X}(r)$ the set of all Delone sets
of $\RR^d$ for which the minimum separation between distinct
points is at least $r$. We assume in the rest of this section
that $r>0$ has been fixed.
 
We define a Hausdorff topology on ${\cal X}(r)$ as follows:
two sets
$S, S' \in {\cal X}(r)$ are ``close'' if for some large compact set 
$K \subset \RR^d$ and some small $\epsilon$  we have
\be \label{topology}
(v + S) \cap K = S' \cap K
\ee
 for some $v \in \RR^r$ with
$|v| < \epsilon$. More precisely we define a uniformity
${\cal U}$ on ${\cal X}(r)$ using as the sets $U(K,\epsilon)$ of uniformity
the set of pairs $(S,S')$
satisfying (\ref{topology}).

\begin{theorem} \cite{RW, boris, martin3}
With respect to this topology ${\cal X}(r)$ is a complete
Hausdorff space.
\end{theorem}

Let $S\in {\cal X}(r)$. Then $\RR^d$ acts on ${\rm LI}(S)$ by
translation and in particular the entire orbit
$[S]$ of $S$ lies in ${\rm LI}(S)$. 
This action is continuous and hence extends also to
an action on the closure $\overline{{\rm LI}(S)}$
of ${\rm LI}(S)$. The relationship between orbits and LI classes
can be summed up by
\be
 S \in [S] \subset {\rm LI}(S) \subset \overline{[S]}
= \overline{{\rm LI}(S)}.
\ee
The second inclusion follows easily from the definitions.
The inclusions may, according to the situation, be
strict or actual equalities.
For a lattice there is only one orbit in its LI class. For
general model sets the situation is very different, as we shall see. 

Recall that a  Delone set $S$ is said to be of 
{\em finite local complexity}
if the closure of $S-S$ is discrete. 
Finite local complexity is a property that is inherited 
by whole LI classes. 

\begin{theorem}  An LI class is pre-compact 
(i.e its completion is compact) iff it has 
finite local complexity.
\end{theorem}

Thus, if  $S$ is a Delone set of finite local complexity 
we obtain a dynamical system
${\cal D}(S)$:
\be \label{dynamicalsystemLI}
\RR^d \times \overline{[S]} \longrightarrow \overline{[S]}.
\ee
In the sequel we will use the symbols like
 ${\cal D}(S)$
to denote both the dynamical system itself and the corresponding
defining space $\overline{[S]}$. 

\begin{theorem} Let $S$ be a Delone set of finite local
complexity. The following are equivalent:
\begin{itemize}
\item[(i)] $S$ is repetitive;
\item[(ii)] $[S] = {\rm LI}(S)$ is closed;
\item[(iii)] The dynamical system ${\cal D}(S)$ is minimal.
\end{itemize} 
\end{theorem}

We recall that minimal means that every $\RR^d$ orbit is dense.

Since generic model sets are repetitive, this leads to
a very nice result:

\begin{theorem}
Let $\Lambda$ be a generic model set.
Then its LI class ${\rm LI}(\Lambda)$ is a compact
Hausdorff space and under the action of translation under
$\RR^d$ it becomes a minimal dynamical system, ${\cal D}(\Lambda)$. 
\end{theorem}

This is the first of the dynamical systems that we wish to 
consider. Its rather abstract form is better understood
by relating it to a more accessible dynamical system.

To this end, let 
$\Lambda = \Lambda(W) = \{ x \in L \mid x^* \in W \}$
be a model set. 
Each element $(u,v)$ of the group $\RR^d \times G$ can be used to 
form a new model set 
\be  
\Lambda(W,u,v) \; := \; u + \{x \in L \mid x^* \in -v+ W \}. 
\ee
If $(u,v) \in {\tilde L} $ then $v = u^*$ and we can rewrite this as 
$\{u + x \in L \mid (u+x)^* \in W \}$, which is just $\Lambda$
again. Thus we get a whole family of model sets parametrized
by $\TT := (\RR^d \times G)/{\tilde L}$ with $\RR^d \times G$ 
acting on it. This is the second dynamical system. Its points
correspond to the model sets $\Lambda(W,u,v)$.
This is the so-called {\em torus parametrization} introduced
by Baake et al. in \cite{BHP}.
We use the same terminology in the more general context here, although
in general $\TT$ is not a torus!

The action of $\RR^d$ on $\TT$,
$ (x , y + {\tilde L}) \mapsto  x+ y + {\tilde L}$,
is a faithful transcription of the operation of translation in physical space, 
so the orbits of $\RR^d$ on $\TT$ correspond to model sets that differ
only by translation. The action of $G$ on $\TT$
corresponds to translating the window around.

\begin{theorem} 
Let $\Lambda$ be a generic model set. Then
\be  \label{dynamicalsysTorus}
\RR^d \times \TT \longrightarrow \TT
\ee
is a minimal uniquely ergodic dynamical system ${\cal D}_{\rm tor}$. 
The unique invariant probability measure is normalized Haar measure.
The set of points of ${\cal D}_{\rm tor}$ corresponding to generic 
model sets is dense and indeed the set of points corresponding
to the non-generic model sets is of the first category.
\end{theorem}

It is noteworthy that
this dynamical system is independent of $W$ but the actual
parametrization of model sets is clearly dependent on it.

So now given a generic model set $\Lambda$, there are two dynamical
systems for the group $\RR^d$, one ${\cal D}(\Lambda)$
coming from the closure of the 
orbit of $\Lambda$
under action of $\RR^d$ and another ${\cal D}_{\rm tor}$
coming from the torus parametrization. Not surprisingly they
are related, but rather surprisingly this relation is somewhat
subtle. All 
the elements of ${\cal D}(\Lambda)$ are, by
definition, in the same LI class. The
same is not the case for the model sets parametrized
by ${\cal D}_{\rm tor}$. Indeed, $\Lambda$
is generic, but translating the window around is bound to 
produce model sets that are not generic. These 
non-generic model sets are not locally isomorphic to the
regular ones, because they have certain special local configurations
of points that are related to the boundaries of their
windows.

\begin{theorem} {\rm \cite{martin3}} Let $\Lambda$ be a generic
model set.
Then there is a continuous surjective
mapping
\be
\beta: {\cal D}(\Lambda) \longrightarrow {\cal D}_{\rm tor}
\ee
which is $\RR^d$-equivariant and which maps $\Lambda$ onto
the point $0$ of the torus. Furthermore, for each of the points of
${\cal D}_{\rm tor}$ which parametrize generic model sets,
the preimage in ${\cal D}(\Lambda)$ consists of a unique point.
\end{theorem}

This mapping comes about as follows: Let $\Lambda' \in {\cal D}(\Lambda)$.
First suppose that $\Lambda' \subset L$. Then it is not hard
to see that $\bigcap_{x\in \Lambda'} (W - x^*)$ is a single point,
call it $b(\Lambda')$. Furthermore, for all $u \in L$, $b(\Lambda' - u)
= b(\Lambda') +u^*$. Now for arbitrary $\Lambda' \in {\cal D}(\Lambda)$,
we can always find $v \in \RR^d$ with $\Lambda' -v \subset L$.
This $v$ is nothing like unique but it follows from what we have just said
that the pair $\beta(\Lambda') := (v, b(\Lambda' -v))$ is unique 
${\rm mod} ~{\tilde L}$, and this is the mapping that we require.

Using these facts it can be established that 

\begin{theorem}
{\rm \cite{martin3}}\label{transferOfStructure}
Assume that $\Lambda$ is a regular and generic
model set. Then
${\cal D}(\Lambda)$ is uniquely ergodic and furthermore
$L^2({\cal D}(\Lambda))$ and $L^2({\cal D}_{\rm tor})$
are isometrically isomorphic as $\RR^d$-spaces.
\end{theorem}

The importance of this is that it shows that from 
the spectrum of ${\cal D}_{\rm tor}$ being discrete, which
it surely is since $\TT$ a compact abelian group, it follows that the 
spectrum of ${\cal D}(\Lambda)$ is also discrete. It is from this
that the pure point diffractivity of $\Lambda$ can be deduced.
In the final section we briefly describe how this happens.

\section{Diffraction}

The theoretical framework for the discussion of diffraction
has been very well described in several places.
The two papers of A.~Hof \cite{bert1, bert2} are standards
and there are also good descriptions in \cite{GK, BMP}. Here we just
quickly formulate the definitions.

Let $\Lambda$ be a regular model set
and define the (tempered) distribution
\be \label{pointmass}
 \delta_\Lambda := \sum_{x \in \Lambda} \delta_x\; ,
\ee
where $\delta_x$ is the Dirac measure at $x$. For each
$s>0$ we calculate the {\em auto-correlation} of $\delta_\Lambda$ restricted
to the ball of radius $s$:
\be
\delta_{\Lambda \cap B_s(0)} * \tilde{\delta}_{\Lambda \cap B_s(0)} 
\; = \; 
\sum_{x,y \in \Lambda \cap B_s(0)}\delta_{x-y} \; ,
\ee
where, as usual, the over-tilde indicates changing the 
sign of the argument. The limit
as $s$ goes to infinity of the volume-averaged auto-correlation
of this measure, which exists for model sets, is the 
{\em auto-correlation measure} of $\Lambda$
(its so-called {\em Patterson function}):
\be
   \gamma \; = \; \lim_{s \to \infty } \frac{1}{\vol(B_s(0))}
   \sum_{x,y \in \Lambda_s} \delta_{x-y}~.  \nonumber
\ee
This limit, taken in the vague topology, converges to a tempered
distribution (i.e. this limit exists when taken against rapidly
decreasing test functions).
Its Fourier transform is a positive measure $\hat{\gamma}$ 
(a result of Bochner's theorem applied to the positive
definite distribution $\gamma$) which
is the {\em diffraction pattern} of $\Lambda$. The measure decomposes
into a point part and a continuous part.
The point part of this measure is the {\em Bragg spectrum} of $\Lambda$. 
The model set has {\em pure point spectrum}
if the continuous part is the trivial $0$-measure.

The complexity of the definition makes it hard to discover
the nature of the diffraction pattern, in particular whether
or not we have pure-point diffraction or not. One approach
has been to use the ergodic theory outlined above, and indeed
it is able to give the main result:

\begin{theorem}\label{main}
{\rm  \cite{martin3}} Any regular model set has pure point spectrum.
Furthermore this spectrum is supported on the projection
into Fourier space on the physical side of the dual
of the compact group $\TT$ {\rm (\ref{dualProj})}, i.e.~ it has the form
\be
{\hat\gamma} \; = \; \sum_{k\in {\hat \TT }}w(k) \delta_{{\hat \pi}_1(k)} 
\ee
\end{theorem}

The proof of this is based on an idea of Dworkin \cite{dworkin}.
The argument is spelled out in \cite{bert2} and we repeat it here
since 
otherwise it is difficult to see the connection between dynamical
systems and diffraction. 

We can assume that $\Lambda$ is generic since translation of
the window does not alter the qualitative nature of the diffraction.
The next step is to replace $\delta_\Lambda$ by a smooth 
approximation to it. To this end, let $b:\RR^d \longrightarrow
\RR_{\ge 0}$ be a smooth function whose support is contained
in the ball $B_r(0)$ of radius $r$, where $B_{2r}(0) \cap
(\Lambda - \Lambda) = \{0\}$. Define a function 
$\psi: {\cal D}(\Lambda) \longrightarrow \RR$ by
\be
\psi(\Lambda') \; = \; \int_{\RR^d}b(-u)\delta_{\Lambda'}(u) du,
\ee
which is continuous on ${\cal D}(\Lambda)$. The action of
$\RR^d$ on $[\Lambda]$ gives rise to a corresponding
action $x \mapsto T_x$ on the space of functions on the orbit $[\Lambda]$
of $\Lambda$ under translation. 

For each $x \in \RR^d$ we have
\be
T_x \, \psi(\Lambda) = \int_{\RR^d}b(-u)\delta_{-x + \Lambda}(u) du
= \int_{\RR^d}b(-u)\delta_{\Lambda}(x+u) du = b*\delta_\Lambda(x) \; ,
\ee
which shows that the function $\sigma^{(b)}\; : \RR^d \longrightarrow
\CC$ defined by 
$x \mapsto T_x(\psi)(\Lambda)$
is obtained by centering a copy of $b$ at each point of $\Lambda$.

Now consider the autocorrelation of $\sigma^{(b)}$:
\be
\begin{array}{ccc}
\gamma^{(b)}(x) &=& \lim_{s \to \infty } \frac{1}{\vol(B_s(0))}
\int_{B_s(0)}\overline{T_{x+y}(\psi)(\Lambda)} T_y(\psi)(\Lambda) dy \\
\\
&=& 
\int_{{\cal D}(\Lambda)}\overline {T_x(\psi}) \psi \, d\mu 
\; = \; (T_x\psi,\psi).
\end{array}
\ee
The main point here is the use of the {\em Birkhoff ergodic theorem}
and the ergodicity of the action of $\RR^d$ on
${\cal D}(\Lambda)$
to replace the integral over $\RR^d$ by an integral over ${\cal D}(\Lambda)$.
Note that the {\em uniqueness} of ergodicity and the continuity of $\psi$
is needed here to guarantee the statement {\em for all} $x$ rather
than a.e. (\cite{Furstenberg}, Sec. 3.2).

In view of Theorem \ref{transferOfStructure}, we have a Fourier expansion
of $\psi$ in terms of the eigenfunctions for the action of $\RR^d$:
$\psi \;= \; \sum a_\lambda \phi_\lambda$, where 
$\phi_\lambda$ is the eigenfunction for the character $x \mapsto 
e^{2\pi i \lambda.x}$ on $\RR^d$.  Thus
$( T_x\psi,\psi) = \sum |a_\lambda|^2 e^{2\pi i \lambda.x}$ and 
taking  Fourier transforms we have 
\be
{\hat \gamma^{(b)}} = \sum |a_\lambda|^2 \delta_\lambda
\ee 
which is a pure point measure on $\RR^d$.

On the other hand we know that $\sigma^{(b)} = b*\delta_\Lambda$
whose autocorrelation can be calculated directly as $b*{\tilde b} *\gamma$
and so ${\hat \gamma^{(b)}} = |\hat b|^2 {\hat \gamma}$. Finally,
taking a sequence of bump functions $\{b\}$ converging in the vague
topology to $\delta_0$ we obtain the required pure-point nature of
the diffraction pattern.

Theorem \ref{main} is qualitative in nature. The quantitative
counterpart is this:

\begin{theorem} \cite{YM2}
Let $k \in {\hat \TT}$ and let $\chi$ denote the
characteristic (or indicator) function of $W$. Then
$w(k) = |{\hat \chi}(- {\hat {\pi}_2}(k))/{\rm vol}(W)|^2$.
\end{theorem} 

There are a number of variations on this theme that are
worthwhile mentioning. First we may imagine replacing the
simple sum (\ref{pointmass}) by a weighted sum
\be \label{weightedpointmass}
 \delta_\Lambda^{\omega} := \sum_{x \in \Lambda} \omega(x)~ \delta_x\; ,
\ee
where $\omega:\Lambda \longrightarrow \CC$ is some function.
\begin{theorem} \label{weightedDiff}
If $\Lambda$ is a regular model set and if $\omega \in {\cal C}(\Lambda)$
(see (\ref{pseudoContinuous})) then the weighted point distribution 
is pure point diffractive.
\end{theorem}

Next we consider the case  that our points of $\Lambda$ are
considered stochastically:
\be \label{stochasticpointmass}
 \delta_{\rm stochastic} := \sum_{x \in \Lambda} \eta(x)~ \delta_x\; ,
\ee
where the $\eta(x)$ form a collection of independent identically
 distributed
random variables that take the  values $1, 0 $ (indicating
occupancy or not of the respective model set sites) with the 
probability of occupancy being $p$:

\begin{theorem}\cite{stochastic}
Let $\Lambda$ be a regular model set and let $\eta$ be as above
with the mean and second moment equal to $m_1$ and $m_2$ 
respectively.
Then the autocorrelation of $\Lambda$ and that of its stochastic version 
are, with probability one, related by
\begin{equation}
        \gamma_{\rm stochastic} \; = \; 
          (m_1)^2 \, \gamma + d~(m_2 - (m_1)^2) \delta_0 \, 
\end{equation}
with Fourier transforms
\begin{equation}
       \hat{\gamma}_{\rm stochastic} \; = \; 
          (m_1)^2 \, \hat{\gamma} + d~(m_2 - (m_1)^2) \, , 
\end{equation}
where $d$ is the density per unit volume of the $\Lambda$.
\end{theorem}

Thus the pure point nature of $\Lambda$ is affected
by at most the addition of a constant continuous background.
More on this stochastic approach may be found in
\cite{stochastic}.

Yet a different variation is to allow the points
of $\Lambda$ to be moved in some regular way.

\begin{theorem}\cite{bert2}
Let $\Lambda$ be a regular model set and let
$f:x \mapsto (f_1(x), \dots, f_d(x))$ be some
mapping of $\Lambda$ into $\RR^d$, where each 
$f_i \in {\cal C}(\Lambda)$ (see (\ref{pseudoContinuous})).
Then the set $\Lambda_f := \{ x + f(x) \mid x \in \Lambda \}$
is pure point diffractive. 
\end{theorem} 

These variations can be combined in the obvious ways.

The problems of determining which point sets are pure
point diffractive is a fascinating and challenging one
which is still wide open.
The examples above show how much model sets can be
modified without serious damage to their diffractive
properties. Obviously adding or removing points
whose average density is $0$ also does not
alter the diffraction. But there are point sets
that are even more remote that are diffractive.
One example is the set of visible points of a lattice.
Given a lattice $L$ in $\RR^d$ its {\em visible} points
are those points $x \in L$, $x \ne 0$,
 satisfying $\QQ x \cap L = \ZZ x$. What is interesting
about the visible points is that they do {\em not}
form a Delone set (they are not relatively dense in $\RR^d$).
In fact, for each $r>0$, the set of holes of radius exceeding
$r$ has positive density. 
However,  
\begin{theorem}\cite{BMP} 
The set of visible points of any lattice of rank at least $2$
is pure point diffractive.
\end{theorem}

\subsection{Comment}
In this paper we have considered point sets in $\RR^d$
that are constructed through the method of projection
from an embedding group and a lattice. We have assumed
that the embedding group is of the form $\RR^d \times G$
where $G$ is a locally compact abelian group, which, in view
of Theorem \ref{MeyerTheorem}, is fairly natural. However,
it is possible to study model sets in the situation
where the `physical space' has been replaced
by an arbitrary locally compact abelian group without
losing many of the most interesting properties. In particular
the diffraction results of Theorem \ref{main} have formulations
in this generality \cite{martin3}.

\section*{Acknowledgment} It is a pleasure to thank Martin Schlottmann
for his edifying insights into this material.

\end{document}